\documentclass{article}  
\usepackage{makeidx}
\usepackage{latexsym,amsmath,amssymb,amsthm}
\usepackage{graphicx}
\title{Generalized Analytic Automorphic Forms for some Arithmetic Congruence subgroups of the Vahlen group on the $n$-Dimensional Hyperbolic Space}
\author{Rolf S\"oren Krau\ss{}har\thanks{Department of
Mathematical Analysis, Ghent University, Building
S-22, Galglaan 2, B-9000 Ghent, Belgium. E-mail:
{\tt krauss@cage.ugent.be}}} \date{9 January, 2004}
\begin{document} 
\maketitle 
\begin{abstract} 
This paper deals with a new analytic type of vector- and Clifford algebra valued automorphic forms in one and two  vector variables. For hypercomplex generalizations of the classical modular group and their arithmetic congruence  subgroups Eisenstein- and Poincar\'e type series that are annihilated by Dirac operators, and more generally, by iterated Dirac operators on the upper half-space of $\mathbb{R}^n$ are discussed. In particular we introduce  (poly-)monogenic  modular forms on hypercomplex generalizations of the classical theta group.           
\end{abstract}
\medskip\noindent

\medskip
 
\newtheorem{definition}{Definition}
\newtheorem{lemma}{Lemma}
\newtheorem{proposition}{Proposition}
\newtheorem{theorem}{Theorem}

\def\symdiff{\triangle} \def\sn{\mathop{\rm
sn}\nolimits} \def\cn{\mathop{\rm
cn}\nolimits} \def\cd{\mathop{\rm cd}\nolimits}
\def\dn{\mathop{\rm dn}\nolimits} \def\BR{{\bf
R}} \def\BC{{\bf C}} \def\BH{{\bf H}} \def\N{{\bf
N}} \def\Z{{\bf Z}} \def\Sc{\mathop{\rm
Sc}\nolimits} \def\Vec{\mathop{\rm Vec}\nolimits}

\section{Introduction}

Automorphic forms are, roughly speaking, functions that show a (quasi-) invariance behavior under the action of a discrete group. For a number of reasons one is in particular interested in those that are furthermore endowed with nice analytic properties. As classical example serve the elliptic modular forms on the modular group $\Gamma:=SL(2,\mathbb{Z})$ which are holomorphic functions on the upper half-plane $H^{+}(\mathbb{C}) = \{z = x+iy \in \mathbb{C}\;|\; y > 0\}$ with a certain boundary condition at $i \infty$ that satisfy 
$$
f(z) = (cz+d)^{-k} f\left(\frac{az+b}{cz+d} \right)  
$$ 
for a positive even integer $k \in 2 \mathbb{N}$ for all $z \in H^{+}(\mathbb{C})$ and all $M = \left(\begin{array}{cc} a & b \\ c & d \end{array} \right) \in \Gamma$.  
The classical examples are the Eisenstein series  
$$
G_k(z)= \sum\limits_{(c,d) \in \mathbb{Z} \times \mathbb{Z} \backslash \{(0,0)\}} (cz+d)^{-k} \quad \quad k \equiv 0(mod\;2),\;\;\;k \ge 4
$$
which first appeared systematically in works of G.~Eisenstein from 1847 and in lectures of K.~Weierstra{\ss} from 1863.
The systematic development of the classical theory of elliptic modular forms on $SL(2,\mathbb{Z})$ is basically due to  H.~Poincar\'e, F.~Klein and R.~Fricke and has been established in the 1890s. In the very beginning one also started to already consider more generally holomorphic modular forms in one complex variable on arithmetic subgroups of the modular group of a finite index. The classical examples are the principal congruence subgroups of level $N \in \mathbb{N}$, 
$$
\Gamma[N] = \Big\{ \left( \begin{array}{cc} a & b \\ c & d \end{array} \right),\; a-1,b,c,d-1 \equiv 0\;(mod\;N) \Big\} 
$$ 
and those subgroups of $\Gamma$ which contain $\Gamma[N]$ as a subgroup. All these groups are called congruence subgroups of $\Gamma$. These include the famous theta group; the complex extension of Euler's historical theta series, i.e.   
$$
\theta(z) = \sum\limits_{n \in \mathbb{Z}} e^{-n^2 \pi z}
$$ 
from 1748 is the prototype of automorphic form on this group involving however a slighlty more general automorphy factor. For detailed information on the classical theory we refer for instance to \cite{schoeneberg} and \cite{KK} and the references therein  in which one can recover a lot of the historical information.\\[0.2cm]  
Extensions of this theory have been developed in several directions. O.~Blumenthal was one of the first to consider in 1904 complex-valued modular forms in several complex variables on Cartesian products of the upper half-plane that show an invariance behavior under hyperabelian groups, direct products of the modular group. In the 1930s C.L.~Siegel and his school started to study holomorphic modular forms in several complex variables for the Siegel modular group which acts on the Siegel half-space of $\mathbb{C}^n$. In 1949 H.~Maa{\ss} introduced in~\cite{Maa49} scalar-valued non-analytic automorphic forms which are eigenfunctions of the hyperbolic Laplace-Beltrami operator on the higher dimensional hyperbolic space and that show an invariance behavior under arithmetic subgroups of the Vahlen group. His work had a remarkable impact. Afterwards many authors started to work in this line of investigation, as for instance A.~Krieg~\cite{Kri2} (1990) and J.~Elstrodt, F.~Grunewald and J.~Mennicke~\cite{EGM90}, V.~Gritsenko~\cite{Grit} among many others.\\[0.2cm] 
In this paper we deal with a new analytic type of vector- and Clifford valued automorphic forms. For a number of important arithmetic congruence groups of the Vahlen group Eisenstein- and Poincar\'e type series are constructed that are null-solutions to the Euclidean Dirac equation and to iterated Euclidean Dirac equations. 

This paper gives a deeper survey as well as an extension of our recent research results presented  in our forthcoming book~\cite{KraHabil} on this topic. In this work we finally come to treat all the congruence groups of the special hypercomplex modular group, that is the group generated by the inversion and the translation matrices inducing the transformations $x \mapsto x+ e_1,\ldots,x\mapsto x+ e_{n-1}$, in the context of Clifford algebra valued function classes in kernels of (iterated) Euclidean Dirac operators. This includes hypercomplex generalizations of the classical theta group. 

\section{Analyticity in Hypercomplex Spaces}

Hypercomplex numbers are, roughly speaking, generalizations of the complex numbers with several imaginary units. While a complex number represents a two dimensional vector in the plane, a hypercomplex number represents a vector in $\mathbb{R}^n$. The advantage of describing two dimensional vectors in terms of complex numbers consists of having an additional multiplication operation available. In dimensions $n > 2$ one can embed $\mathbb{R}^n$ into the so-called {\it real Clifford algebra} over ${\mathbb{R}}^n$ which then allows us to also endow $\mathbb{R}^n$ with a multiplicative structure. In this section we briefly recall some basic notions on Clifford algebras and their function theory. For details we refer the interested reader for example to~\cite{BDS}. 

\subsection{Clifford algebras}

Throughout this article $\{e_1,e_2,\ldots,e_n\}$ stands for the
canonical basis vectors of the Euclidean space
${\mathbb{R}}^n$, equipped with the quadratic form $Q(x) = x_1^2+\cdots+x_n^2$.\\ Here, and in all that follows, the letter $x$ stands for a vector from ${\mathbb{R}}^n$, i.e. $x = x_1 e_1 + \cdots + x_n e_n$.\\  
The attached {\it real Clifford algebra} 
$Cl_n$ of the quadratic space $({\mathbb{R}}^n,Q)$ is then the free algebra that is 
generated by ${\mathbb{R}}^n$ modulo the relation 
\begin{equation}
\label{clifford}
x^2 = - Q(x) e_0,
\end{equation}
where $e_0$ denotes the neutral multiplicative element of the Clifford
algebra $Cl_n$.\\
From (\ref{clifford}) one derives directly the following multiplication 
rules for the elements of the canonical basis from ${\mathbb{R}}^n$: 
\begin{equation} 
e_i e_j + e_j e_i = - 2 \delta_{ij} e_0,\quad i,j=1,\cdots,n. 
\end{equation}
Here, $\delta_{ij}$ denotes the Kronecker symbol. 
A vector space basis for the whole Clifford algebra $Cl_n$ is given by the
set 
$$
\{e_A : A \subseteq \{1,\cdots,n\}\} \quad \mbox{with} \quad e_A = e_{l_1} e_{l_2} \cdots e_{l_r},
$$ 
where 
$$
1 \le l_1 < \cdots < l_r \le n,\;\; e_{\emptyset} = e_0 = 1.
$$
Every  $a \in Cl_n$ can be written in the form $ a = \sum\limits_{A \in P(\{1,\ldots,n\})} a_A e_A$ with
$a_A \in {\mathbb{R}}$. Two examples of real Clifford algebras are the complex number 
field ${\mathbb{C}}$ and the Hamiltonian skew field ${\mathbb{H}}$. In this sense one can regard Clifford algebras as higher dimensional associative generalizations of the complex numbers.\\[0.2cm]
We further consider the {\it conjugation} anti-automorphism in the Clifford algebra $Cl_{n}$ which is defined by
$\overline{a} = \sum_A a_A \overline{e}_A$,
where $\overline{e}_A = \overline{e}_{l_r}
\overline{e}_{l_{r-1}} \cdots \overline{e}_{l_1}$
and $\overline{e}_j = - e_j$ for $j=1,\cdots,n,\;
\overline{e}_0 = e_0 = 1$.\\[0.2cm]
The {\it reversion} anti-automorphism is defined by $a^{*} := \sum_A (-1)^{|A| \frac{|A|-1}{2}}   a_A e_A$, 
where $|A|$ denotes the cardinality of the set $A$. Similar to the conjugation automorphism, the reversion invertes the order of the factors in a product of Clifford numbers, i.e. $(ab)^{*}= b^{*}a^{*}$. However, it does not produce a minus sign when applied to a vector, that is $e_i^{*} = e_i$.\\[0.2cm]
The subspace ${\cal{A}}_{n+1} := \mbox{span}_{\mathbb{R}}
\{1,e_1,\cdots,e_n\} = {\mathbb{R}} \oplus {\mathbb{R}}^n
\subset Cl_{n}$  is often called the space
of paravectors $z = x_0 + x_1 e_1 + x_2 e_2 +
\cdots + x_n e_n$. Each paravector $z \in {\cal{A}}_{n+1}$ satisfies $z^{*} = z$ and for each vector $x \in \mathbb{R}^n$ holds $\overline{x} = -x$.    
One can identify the spaces ${\cal{A}}_{n+1}$ with ${\mathbb{R}}^{n+1}$; we prefer to  work in this article with vectors instead of paravectors for the sake of simplicity. For our needs the treatment with vectors has the advantage that we can exploit the following relationship $x^{*} = - \overline{x} = x$ for all $x \in \mathbb{R}^n$ when we work with vectors.    

A (pseudo-) scalar product between two Clifford numbers
$a,b \in Cl_n$ is defined by 
$$
\langle a,b \rangle := \Sc(a\overline{b})
$$   
and the associated Clifford (pseudo-) norm of an arbitrary  $a = \sum\limits_A a_A e_A$
is 
$$
\|a\| = ( \sum\limits_A |a_A|^2)^{1/2}.
$$
Each vector $x \in \mathbb{R}^{n} \backslash \{0\}$ has an inverse element
in ${\mathbb{R}}^{n}$ given by $x^{-1} = -\frac{x}{\|x\|^2}$. 

\subsection{Generalized analyticity in hypercomplex spaces}
There are several ways to generalize complex analyticity to higher dimensions. One approach is offered by 
Clifford analysis which considers Clifford algebra valued functions that solve Cauchy-Riemann or Dirac type equations in regions of higher dimensional vector spaces, as for instance in ${\mathbb{R}}^{n}$ in the simplest case, or, more generally on some manifolds.  In this paper we restrict to considering open subsets in the Euclidean space ${\mathbb{R}}^{n}$.\\  
In the Euclidean space ${\mathbb{R}}^{n}$ the attached  
{\it Dirac operator} is simply given by   
$$
D_x := \sum\limits_{i=1}^n
e_i\frac{\partial }{\partial x_i}.
$$ 
When it is clear to which variable the operator is applied, then we simply write $D$ for $D_x$. 
In the case $n=2$, the operator $-D e_2$ can be identified with he classical complex Cauchy-Riemann operator.\\[0.2cm]  
Following~\cite{BDS} and other references, a $Cl_n$-valued real differentiable function $f$ that is defined in an open subset $U \subseteq {\mathbb{R}}^{n}$ is called {\it left (right) monogenic} 
at a point $z \in U$ if $D_xf(x) = 0$ or $ f(x)D_x = 0$, respectively.\\[0.2cm] 
The notion of left (right)
monogenicity in $\mathbb{R}^{n}$ provides hence  
a higher dimensional generalization of the concept of complex
analyticity in the sense
of the Riemann approach. Indeed, a number of 
classical theorems from complex analysis could
be generalized to higher dimensions by this
approach, including in particular an analogue of Cauchy's theorem and a Cauchy integral formula. See for instance~\cite{BDS}.\\[0.2cm] Clifford analysis may also be regarded as a special function theory within harmonic analysis: The Dirac operator $D$ factorizes the Euclidean Laplacian $\Delta = \sum_{i=1}^n \frac{\partial^2 }{\partial x_j^2}$ viz $D^2 = - \Delta$. Each real component of a left (right) monogenic function is hence Euclidean harmonic.  In this paper we also deal more generally with {\it left (right) $s$-monogenic} functions, or {\it polymonogenic} functions for short,  which are sufficiently smooth $Cl_n$-valued functions that are annihilated by an iterated Dirac operator, that means that satisfy $D^s f = 0$ resp. $f D^s = 0$ for a positive integer $s \in \mathbb{N}$. For even $s$ the operator $D^s$ coincides up to a minus sign with the iterated Euclidean Laplacian $\Delta^{s/2}$. Following for instance \cite{Ry93,Ry2000}, the fundamental solution to $D^s$ is given in the cases where $s < n$ by 
\begin{equation}
q_{\bf 0}^{(s)}(x) := \left\{ \begin{array}{cc} \frac{x}{\|x\|^{n+1-s}}  & s \mbox{ odd integer with } s < n,\\
\frac{1}{\|x\|^{n-s}}  & s \mbox{ even integer with } s < n, 
\end{array} \right.
\end{equation} 
The functions serve as Green kernels and give rise to Green type formulas for $s$-monogenic function, as described in detail for instance in~\cite{Ry2000}. 
The fundamental solution to the linear Dirac operator reads thus $q_{\bf 0}(x) = - \frac{x}{\|x\|^n}$ and generalizes  the classical Cauchy kernel function to higher dimensions. 

For the sake of clarity we use multi-index notation throughout this paper. In this sense, we write  ${\bf m} = (m_1,\cdots,m_n)$ for a multi-index from ${\mathbb{N}}_0^n$; the expression $|{\bf m}|= m_1+\cdots+m_n$ abbreviates its length. We further write $\tau(j)$ for the particular index where $m_i = \delta_{ij}$ for $i=1,2,\ldots,n$, $\delta_{ij}$ standing for the Kronecker symbol. For a vector $x\in {\mathbb{R}}^n$, we write $x^{\bf m} = x_1^{m_1} \cdot \cdots \cdot x_n^{m_n}$. Consequently we write 
\begin{equation}
q^{(s)}_{\bf m}(x) := \frac{\partial^{|{\bf m}|}}{\partial x^{\bf m}} q^{(s)}_{\bf 0}(x) 
\end{equation}
for the partial derivatives of $q^{(s)}_{\bf 0}$ with respect to differentiations in the $x_1,\cdots,x_n$ directions.
\section{The Vahlen group and some basic geometric and arithmetic subgroups}

As very well-known, M\"obius transformations in the plane can be represented by $2 \times 2$-matrices from $GL(2,{\mathbb{C}})$. In 1902 K.~Th.~Vahlen~\cite{Vahlen} discovered that one can describe M\"obius transformation in ${\mathbb{R}}^{n}$ in a similar way in terms of a matrix group which consists of special $2 \times 2$ Clifford matrices. For convenience, let us recall here its definition:
\begin{definition}\label{GV}(General Vahlen group)\\
A matrix $M = \left(\begin{array}{cc} a & b \\ c & d \end{array} \right)$ with $Cl_n$-valued entries $a,b,c,d$ belong to the general Vahlen group $GV({\mathbb{R}}^{n})$ if $a,b,c,d$ can all be written as products of vectors from $\mathbb{R}^{n}$ and if  
\begin{eqnarray}
ad^{*} - bc^{*} & \in & \mathbb{R} \backslash\{0\} \\
a^{-1} b, c^{-1} d & \in & \mathbb{R}^{n},\;\mbox{if}\;c\neq 0\;\mbox{or}\; a \neq 0,\;\mbox{respectively}.
\end{eqnarray}
\end{definition}
Indeed, if $T: {\mathbb{R}}^{n} \cup \{\infty\} \rightarrow \mathbb{R}^{n} \cup \{\infty\}$ is a M\"obius transformation, then it can be written in the form 
$$
T(x) = (ax+b)(cx+d)^{-1}
$$ 
with coefficients $a,b,c,d$ stemming from a Vahlen matrix from $GV({\mathbb{R}}^{n})$. Since $x^{*} = x$ for all $x \in \mathbb{R}^{n}$ we can also rewrite the expression $(ax+b)(cx+d)^{-1}$ equivalently in the way $(xc^{*}+d^{*})^{-1}(xa^{*}+b^{*})$.\\[0.2cm]  
The geometric properties of Vahlen matrices and their associated M\"obius transformations have been studied extensively, for instance in~\cite{Ahlf}. Every M\"obius transformation is composed by four types of elementary transformations: by translations, rotations, the inversion at the unit sphere and dilatations. Their associated Vahlen matrices have the form 
$$
T_b= \left(\begin{array}{cc} 1 & b \\ 0 & 1 \end{array} \right),\quad R= \left(\begin{array}{cc} u^{*} & 0 \\ 0 & u^{-1} \end{array} \right), \quad J = \left(\begin{array}{cc} 0 & -1 \\ 1 & 0 \end{array} \right), \quad D =  \left(\begin{array}{cc} \alpha & 0 \\ 0 & \delta \end{array} \right),  
$$
where $b \in {\mathbb{R}}^{n}$, $u = u_1 \cdots u_t$ with $u_i \in S^{n} :=\{x \in \mathbb{R}^{n}\;|\; \|x\| = 1\}$, and $\alpha,\delta \in \mathbb{R} \backslash\{0\}$. As shown for example in~\cite{EGM87} and elsewhere, the translation and dilatation matrices together with the inversion generate already the complete general Vahlen group $GV(\mathbb{R}^{n})$. The rotation matrices can be constructed from these three types of matrices.\\[0.2cm]
The set of those Vahlen matrices where the coefficients $a,b,c,d$ satisfy additionally $ad^{*}-bc^{*}=1$ forms a normal subgroup of $GV(\mathbb{R}^{n})$ and will be called the {\it special Vahlen group} $SV({\mathbb{R}}^{n})$. This subgroup can be generated completely by the translation matrices and the inversion, as shown for instance in~\cite{EGM87}. It contains only those dilatation matrices $D$ which have the property that $\delta=1/\alpha$.\\
The subgroup $SV(\mathbb{R}^{n-1})$, consisting of those matrices from $SV(\mathbb{R}^{n})$, where the coefficients $a,b,c,d$ satisfy even $a^{-1} b, c^{-1} d  \in  \mathbb{R}^{n-1}$ if $c\neq 0$ or $a \neq 0$, respectively, acts transitively on the upper half-space 
\begin{equation}
H^{+}(\mathbb{R}^{n}) = \{x \in \mathbb{R}^{n}\;|\;x_n > 0\}.  
\end{equation}
Here we want to focus now on some particular arithmetic subgroups of $SV(\mathbb{R}^{n-1})$ which act discontinuously on the upper half-space.  
A general description of discrete arithmetric subgroups of the special Vahlen group $SV(\mathbb{R}^{n-1})$ is given for instance~in \cite{EGM90}. In this paper we want to consider the {\it special hypercomplex modular group}   
\begin{equation}
\Gamma_p = \Big\langle T_{e_1},\ldots,T_{e_p},J\Big\rangle \quad \quad \mbox{where}\;\; p \in \{1,2,\ldots,n-1\} 
\end{equation}
and their basic associated congruence subgroups. Let $N$ be a positive integer. Then the {\it principal congruence subgroups of level $N$ of $\Gamma_p$} are said to be the groups  
\begin{equation}
{\Gamma_p}[N] = \Big\{ \left( \begin{array}{cc} a & b \\ c & d \end{array} \right) \in \Gamma_p \;\Big|\;\; a-1,b,c,d-1 \in N {\cal{O}}_p \Big\} 
\end{equation}
where ${\cal{O}}_p = \sum_{A \in P(\{1,\ldots,p\})} \mathbb{Z} e_A$ stands for the standard order in $Cl_n$. For all $N \ge 1$ the groups $\Gamma_p[N]$ are normal subgroups of $\Gamma_p$. For $N=1$, the group $\Gamma_p[1]$ coincides with $\Gamma_p$.

A subgroup $\Lambda \subseteq \Gamma_p$ is called a {\it congruence subgroup} of $\Gamma_p$ if there exists an $N \in \mathbb{N}$ such that $\Gamma_p[N] \subset \Lambda \subset \Gamma_p$.  

All congruence groups of the hypercomplex modular group have thus a finite index in $\Gamma_p$. 

Two basic examples are 
\begin{equation}
{\Gamma_p}^{0}[N] = \Big\{ \left( \begin{array}{cc} a & b \\ c & d \end{array} \right) \in \Gamma_p \;\Big|\;\; b \in N {\cal{O}}_p \Big\} 
\end{equation}
\begin{equation}
{\Gamma_p}_{0}[N] = \Big\{ \left( \begin{array}{cc} a & b \\ c & d \end{array} \right) \in \Gamma_p \;\Big|\;\; c \in N {\cal{O}}_p \Big\} 
\end{equation} 
 Note that ${\Gamma_p}[N] \subseteq {\Gamma_p}^{0}[N] \subseteq \Gamma_p$ and  ${\Gamma_p}[N] \subseteq {\Gamma_p}_{0}[N] \subseteq \Gamma_p$ for all $N \in {\mathbb{N}}$.  

A further interesting example is the following group 
\begin{equation}
\label{theta}
{\Gamma_p}_{\theta} = \Gamma_p[2] \cup \Gamma_p[2] J
\end{equation}
which can be regarded as a higher dimensional generalization of the classical {\it theta group}. 

For all $N > 1$, the congruence groups ${\Gamma_p}^{0}[N]$, ${\Gamma_p}_{0}[N]$ are not normal subgroups of $\Gamma_p$, neither so ${\Gamma_p}_{\theta}$. The generalized theta group ${\Gamma_p}_{\theta}$ and the groups ${\Gamma_p}^{0}[2]$, ${\Gamma_p}_{0}[2]$ are conjugated in $\Gamma_p$. 

All these particular groups are indeed discrete subgroups of the Vahlen group $SV(\mathbb{R}^{n-1})$ and act discontinuously on the upper half-space $H^{+}(\mathbb{R}^n)$ by its associated M\"obius transformation. 

In the next section it will be explained how to construct monogenic and polymonogenic automorphic forms for these particular arithmetic groups.  

\section{Polymonogenic modular forms for congruence groups of the special hypercomplex modular group}
In contrast to classical complex analysis, neither the multiplication nor the composition of two (poly-)monogenic functions results into a new monogenic function again, in general. The set of left (right) $s$-monogenic function forms just a right (left) $Cl_n$-module, respectively. However, the Dirac operator and its iterates have a further  important monogenicity preserving property. It is (quasi-) invariant under the action of $SV(\mathbb{R}^n)$. Following for instance \cite{Ry93}: If $\left( \begin{array}{cc} a & b \\ c & d \end{array} \right)$ is a matrix from $SV(\mathbb{R}^n)$ and $f$ a left $s$-monogenic function ($s < n$) in the variable $y:= (ax+b)(cx+d)^{-1}$, then the function $F(x) := q^{(s)}_{\bf 0}(cx+d) f((ax+b)(cx+d)^{-1})$ turns out to be $s$-monogenic in the variable $x$, whenever  $cx+d \neq 0$.\\[0.2cm] 
All the congruence groups that are considered in the previous section act discontinously on the upper half-space. As a consequence of this, the expression $(ax+b)(cx+d)^{-1}$ is for all $x \in H^{+}(\mathbb{R}^n)$ and all matrices $\left( \begin{array}{cc} a & b \\ c & d \end{array} \right)$ from these congruence groups again a  well-defined object from the  half space $H^{+}(\mathbb{R}^n)$.  Therefore, one can strengthen the conformal invariance property in the following form:
\begin{lemma}
Let $s,p < n$ and $\Lambda$ be a subgroup of the hypercomplex modular group $\Gamma_p$. If $f: H^{+}(\mathbb{R}^n) \rightarrow Cl_n$ satisfies the equation $D_x^s f(x) = 0$ at each $x \in H^{+}(\mathbb{R}^n)$, then the function 
$$
F(x) = q^{(s)}_{\bf 0}(cx+d) f((ax+b)(cx+d)^{-1})
$$
satisfies for all $M \in \Lambda$ on the whole half-space equations of the form $D_x^l F(x) = 0$ for all positive integers $l \ge s$. 
\end{lemma}
One gets a similar statement for right $s$-monogenic functions, involving an analogous transformation of the form 
$$
F(x) = f((ax+b)(cx+d)^{-1})  q^{(s)}_{\bf 0}(xc^{*}+d^{*}).
$$ 
In view of $\overline{cx+d} = \pm (cx+d)^{*}$ also a transformation of the form 
$$
F(x) = \overline{q^{(s)}_{\bf 0}(cx+d)^{*}} f((ax+b)(cx+d)^{-1})
$$
remains left $s$-monogenic. A similar statement holds for the right $s$-monogenic case.\\[0.2cm]
Notice that left $s$-monogenic functions form a Clifford right module, that means we can add them without losing $s$- monogenicity.  
As a consequence of Weierstra{\ss}' convergence theorem, series over expressions of the form $q^{(s)}_{\bf 0}(cx+d) f((ax+b)(cx+d)^{-1})$ yield left $s$-monogenic functions on the half-space, provided it converges normally there.\\[0.2cm]    
That is the analytic basis for the construction of $s$-monogenic Eisenstein- and Poincar\'e series. In order to construct examples of left $s$-monogenic automorphic forms on the congruence groups $\Lambda \subseteq \Gamma_p$ it is thus suggestive at the very first to start with a left $s$-monogenic function that is bounded on $H^{+}(\mathbb{R}^n)$ and to sum the expressions  $q^{(s)}_{\bf 0}(cx+d) f((ax+b)(cx+d)^{-1})$ over all matrices from the complete group $\Lambda$. However, in view of having infinitely many matrices with $c=0$, this would lead to a divergent series. To overcome this problem, one possibility is to already start with a bounded $s$-monogenic function $\tilde{f}:H^{+}(\mathbb{R}^n) \rightarrow Cl_n$ that is additionally totally invariant under the subgroup of translation matrices $T(\Lambda)$ that is contained in the congruence group $\Lambda$ under consideration. Then it namely suffices to consider summations over a complete set of representatives of right cosets in $\Lambda$ modulo $T(\Lambda)$ instead of extending the summation over the whole group $\Lambda$. 

In what follows the notation $M: T(\Lambda) \backslash \Lambda$ means that $M$ runs through a system of representatives ${\cal{R}}$ of the right cosets of $\Lambda$ with respect to its translation subgroup $T(\Lambda)$, i.e. 
$$
\cup_{M \in {\cal{R}}} T(\Lambda) M = \Lambda\;\; {\rm and}\;\; T(\Lambda) M \neq T(\Lambda) L \;\;{\rm for}\;\; M,N \in {\cal{R}}\;\; {\rm with}\;\; M \neq L. 
$$
The following theorem provides us with a generalization of \cite{KraHabil}, [Theorem 3.4] to arbitrary congruence subgroups of $\Gamma_p$: 
\begin{theorem} 
\label{firstconstruthm}
Let $n \in {\mathbb{N}}, s \in {\mathbb{N}}, s < n, \; p \in \{1,\ldots,n-1\}$ and $p < n - s -1$.\\ Let $\Lambda \subset \Gamma_p$ be a congruence subgroup and let $T(\Lambda)$ be its subgroup of translation matrices. Further, let $\tilde{f}:H^{+}(\mathbb{R}^n) \rightarrow {Cl}_{n}$ be a bounded and left monogenic function on $H^{+}({\mathbb{R}}^n)$ that is totally invariant under the translation subgroup $T(\Lambda)$. Then
\begin{equation}
\label{constru1}
f(x) = \sum\limits_{M: T(\Lambda) \backslash \Lambda} q^{(s)}_{\bf 0}(cx+d) \tilde{f}((ax+b)(cx+d)^{-1}) \quad \quad x \in H^{+}(\mathbb{R}^n)
\end{equation}
is a Clifford-valued function which is bounded in any compact subset of $H^{+}(\mathbb{R}^n)$. Moreover, it satisfies on the whole upper half-space ${\cal{D}}^l f(x) = 0$ for 
all $l \ge s$ and $f(x) \;=\; q^{(s)}_{\bf 0}(cx+d) f((ax+b)(cx+d)^{-1}) $ for all  $M \in \Lambda$.
\end{theorem}
{\it Sketch of the Proof.}\\ 
Since $\Lambda \subseteq \Gamma_p$, each expression $(ax+b)(cx+d)^{-1} \in H^{+}(\mathbb{R}^n)$ for all $x \in H^{+}(\mathbb{R}^n)$. In view of the boundedness of $\tilde{f}$ on $H^{+}(\mathbb{R}^n)$ we may hence conclude that there is a positive real number $L$ such that     
$$
\|\tilde{f}((ax+b)(cx+d)^{-1})\| \le L 
$$
for all $x \in H^{+}({\mathbb{R}}^n)$ and for all $M = \left( \begin{array}{cc} a & b \\ c & d \end{array} \right) \in \Lambda$.  
Therefore, it suffices to show the normal convergence of 
\begin{equation}
\label{estimate}
\sum\limits_{M: T(\Lambda) \backslash \Lambda} \|c \tau +d\|^{-\alpha}
\end{equation}
on $H^{+}(\mathbb{R}^n)$ for $\alpha > p+1$. 

To this end one considers an arbitrary matrix 
$
\left(
\begin{array}{cc}
* & * \\
c & d
\end{array} \right) 
\in \Lambda$ and for arbitrary  $\varepsilon > 0$ the vertical strip
\begin{equation}
V_{\varepsilon}(H^{+}({\mathbb{R}}^n)) := \Big\{ \tau =(\underline{{\bf x}},x_n) \in H^{+}({\mathbb{R}}^n) \Big|  
\underline{\bf x} \in {\mathbb{R}}^{n-1}:\;
\|\underline{{\bf x}}\| \le \frac{1}{\varepsilon}, x_n > \varepsilon \Big\}.
\end{equation} 
By using a classical compactification argument (for details, see e.g. Proof of [Theorem 3.4] from \cite{KraHabil}) one can show that one can find for every $\varepsilon > 0$  a real $\rho > 0$ so that   
\begin{equation}
\label{rho}
\|c \tau + d\| \ge \rho \|c e_n + d\| \;\;\;\; \forall \tau \in V_{\varepsilon}(H^{+}(\mathbb{R}^n)) \;\;\;\mbox{and}\;\;\;\;\left(
\begin{array}{cc}
* & * \\
c & d
\end{array} \right) 
\in \Lambda
\end{equation}  

Since the index $[\Gamma_p:\Lambda] < \infty$, the series 
\begin{equation}
\label{estimate1}
\sum\limits_{M: T(\Lambda) \backslash \Lambda} \|c e_n + d\|^{-\alpha}
\end{equation}
has the same convergence abscissa as the series 
\begin{equation}
\label{estimate2}
\sum\limits_{M: {\cal{T}}_p \backslash \Gamma_p} \|c e_{n} + d\|^{-\alpha},
\end{equation} 
where ${\cal{T}}_p = \langle T_{e_1},\ldots,T_{e_p}\rangle$ is the maximal translation group contained in the special hypercomplex group $\Gamma_p$. This series in turn has of course has the same convergence abscissa as 
\begin{equation}
\label{estimate3}
\sum\limits_{M: {\cal{T}}_p \backslash \Gamma_p} \|c e_{p+1} + d\|^{-\alpha}
\end{equation}
since the elements $e_{p+1},\ldots,e_n$ play all an equal role satisfying the same calculation rules in relation with the elements $c,d \in span_{\mathbb{R}}\{1,e_1,\ldots,e_p,\ldots,e_{1\ldots p}\}$. According to \cite{EGM90}, the series (\ref{estimate3}) has precisely the convergence abscissa $p+1$. Hence, the series~\ref{constru1} converges normally on $H^{+}(\mathbb{R}^n)$.\\ 
By Weierstra{\ss}' convergence theorem then follows that $f$ satisfies ${\cal{D}}^l f = 0$ for all $l \ge s$ in $H^{+}({\mathbb{R}}^n)$, since $\tilde{f}$ is monogenic in $H^{+}(\mathbb{R}^n)$.\\
To verify that $f$ is an automorphic form with respect to $\Lambda$, consider an arbitrary matrix $L$ from the congruence group $\Lambda$. In the sequel, we write for simplicity $M\langle x \rangle := (ax+b)(cx+d)^{-1}$ where 
$ M = \left(
\begin{array}{cc}
a & b \\
c & d
\end{array} \right)$. 

In view of the relation $q^{(s)}_{\bf 0}(ab) = q^{(s)}_{\bf 0}(b) q^{(s)}_{\bf 0}(a)$, both for $s$ even and odd one hence obtains\\
$
f(L\langle x \rangle)  =  \sum\limits_{M: T(\Lambda) \backslash \Lambda} 
q^{(s)}_{\bf 0}(c_M L \langle x \rangle + d_M) \tilde{f}(M \langle L \langle x \rangle\rangle)
$
\begin{eqnarray*}
& & =  \sum\limits_{M: T(\Lambda) \backslash \Lambda} q^{(s)}_{\bf 0}\Big(c_M(a_L x + b_L)(c_L x + d_L)^{-1} + d_M\Big) \tilde{f}(ML \langle x \rangle)\\
& & =  \sum\limits_{M: T(\Lambda) \backslash \Lambda} q^{(s)}_{\bf 0} \Bigg(
\frac{c_M(a_L x + b_L)(\overline{c_L x + d_L}) + d_M(c_L x + d_L) \overline{(c_L x + d_L)}}{\|c_L x + d_L\|^2}\Bigg) \tilde{f}(ML \langle x  \rangle)\\
& & =  \sum\limits_{M: T(\Lambda) \backslash \Lambda} q^{(s)}_{\bf 0}\Bigg(\frac{\overline{c_L x + d_L}}{\|c_L x + d_L \|^2}\Bigg) \;  q^{(s)}_{\bf 0}\Big((c_M a_L + d_M c_L)x + c_M b_L + d_M d_L\Big) \tilde{f}(ML \langle x \rangle)\\
& & =  \Big[q^{(s)}_{\bf 0}(c_L x + d_L)\Big]^{-1} \sum\limits_{M: T(\Lambda) \backslash \Lambda} q^{(s)}_{\bf 0}(c_{ML}x+d_{ML}) \tilde{f}(ML \langle x \rangle)\\
& & =  \Big[q^{(s)}_{\bf 0}(c_L x + d_L)\Big]^{-1} f(x).
\end{eqnarray*}
The last step follows by a rearrangement argument which may be established in view of the totally invariance of  $\tilde{f}$ is invariant under $T(\Lambda)$.\hspace{1cm} Q.E.D.\\[0.2cm]
The next important step is to analyse for which weights $s$ and for which congruence groups $\Lambda$ this construction yields definitely {\it non-vanishing}  modular forms.\\[0.3cm] 
{\bf Even weights}\\[0.2cm]
The construction principle proposed in Theorem~\ref{firstconstruthm} generates non-trivial Clifford valued examples of automorphic forms on all the congruence groups $\Lambda \subseteq \Gamma_p$ with $p < n-1-s$ within the function class Ker $D^l$ for all positive integers $l \ge s$ whenever $s$ is an even integer with $s < n$.  The simplest examples can be obtained by putting $\tilde{f} \equiv 1$ which leads to the function series 
\begin{equation}
\label{escalar}
{\cal{G}}^{(\Lambda)}_{s} (x) = \sum\limits_{M: T(\Lambda) \backslash \Lambda} \frac{1}{\|c x +d\|^{n-s}}
\end{equation}
The series in (\ref{escalar}) have formally a similar form as the non-analytic Eisenstein series considered in~\cite{EGM90}.\\ 
Following Theorem~\ref{firstconstruthm}, this series converges for all congruence subgroups of $\Gamma_p$ whenever $p < n-1-s$. 
\begin{proposition}
Let $s<n$ be an even integer and suppose that $p < n-1-s$. For each congruence subgroup $\Lambda \subseteq \Gamma_p$ the series  ${\cal{G}}^{(\Lambda)}_{s} (x)$ in (\ref{escalar}) do not vanish identically. 
\end{proposition}
{\it Proof}. 
To show that ${\cal{G}}^{(\Lambda)}_{s} (x) \not\equiv 0$ for all congruence subgroups of $\Gamma_p$ one simply needs to consider the following limit:   
\begin{eqnarray}
\label{lim1}
\lim\limits_{x_n \rightarrow \infty} {\cal{G}}^{(\Lambda)}_{s} (x)  &=&   \sum\limits_{M: T(\Lambda) \backslash \Lambda} \lim\limits_{x_n \rightarrow \infty} \|c x + d\|^{s-n}\nonumber \\
& = & \sum\limits_{M: T(\Lambda) \backslash \Lambda, c \neq 0} \underbrace{\lim\limits_{x_n \rightarrow \infty} \|c x + d\|^{s-n}}_{=0} \nonumber \\
 & + &  \sum\limits_{M: T(\Lambda) \backslash \Lambda,c = 0} \|d\|^{s-k}. 
\end{eqnarray}
We have $\sharp\{M \in T(\Gamma_p) \backslash \Gamma_p, c = 0\} = 2^{p+1}$. Since all congruence subgroups $\Lambda$ in $\Gamma_p$ contain the identity matrix, we have 
\begin{equation}
1 \le \sharp\{M \in T(\Lambda) \backslash \Lambda, c = 0\} \le 2^{p+1} \quad \quad \mbox{for\;all}\; \Lambda.
\end{equation}
The expression in~(\ref{lim1}) is hence different from zero for all congruence subgroups, including of course in particular the principal congruence subgroups $\Gamma_p[N]$ of any level $N \in \mathbb{N}$, the other two basic types of congruence groups ${\Gamma_p}^{0}[N],{\Gamma_p}_{0}[N]$, as well as the generalized theta groups ${\Gamma_p}_{\theta}$. Q.E.D.\\[0.2cm] 
The series ${\cal{G}}^{(\Lambda)}_{s} (x)$ provide hence for all $p < n-1-s$ and all even $s < n$ non-trivial examples of $l$-monogenic automorphic forms for all $l \ge s$ for all congruence subgroups of $\Gamma_p$.  All of them are scalar-valued.\\[0.2cm] 
Inserting for $\tilde{f}$ special variants of the following family of monogenic Eisenstein series~\cite{KraHabil},~[Chapter 2.4] 
$$
G_{\bf m}(x) = \sum\limits_{(\alpha,\omega) \in \mathbb{Z} \times \mathbb{Z}^{n-1} \backslash \{(0,\underline{\bf 0})\}} q^{(1)}_{\bf m}(\alpha x + \omega), \quad \quad  |{\bf m}|\ge 3,\;|{\bf m}| \equiv 1(mod\;2)
$$
leads to a family of non-trivial vector-valued examples in classes of polymonogenic functions for all congruence groups $\Lambda\subseteq \Gamma_p$ with $p < n-1-s$ and all even weights $s$ with $s < n$. Notice that $G_{\bf m}$ is $T(\Lambda)$-invariant, since each translation subgroup $T(\Lambda)$ of an arbitrary congruence subgroup $\Lambda \subseteq \Gamma_p$ is contained in the translation invariance group 
\begin{equation}
{\cal{T}}_{n-1} = \langle T_{e_1},\ldots, T_{e_{n-1}} \rangle 
\end{equation}
of the series $G_{\bf m}(x)$. 
\begin{proposition}
Let $s<n$ be an even positiv integer and suppose that $p < n-1-s$. Let $\Lambda \subseteq \Gamma_p$ be an arbitrary congruence subgroup. 
Then there are multi-indices ${\bf m} \in \mathbb{N}_0^n$ for which the associated  series 
\begin{equation}
E^{\Lambda}_{s,{\bf m}}(x) =  \sum\limits_{M: T(\Lambda) \backslash \Lambda} q^{(s)}_{\bf 0}(cx+d) G_{\bf m}(M \langle x \rangle +e_n)
\end{equation}
do not vanish identically. They represent non-trivial vector-valued modular forms on $\Lambda$.  
\end{proposition}
{\it Proof}. To show that the family of series  $E^{\Lambda}_{s,{\bf m}}(x)$ 
contain indeed non-vanishing candidates, we again can use the previous limit argument. For clarity let us use the notation $x = \underline{{\bf x}} + x_n e_n$ with $\underline{{\bf x}} \in {\mathbb{R}}^{n-1}$.\\[0.2cm] 
For all $x \in H^{+}(\mathbb{R}^n)$ and all matrices $M$ from an arbitrary congruence subgroups of $\Gamma_p$ we have  
$$
M \langle x \rangle + e_n \in \{x \in H^{+}({\mathbb{R}}^n) | x_n \ge 1\}.
$$
We can always find a lattice point  $\omega_0 \in {\mathbb{Z}}^{n-1} = \mathbb{Z} e_1 + \cdots + \mathbb{Z} e_{n-1}$ such   
that $M \langle x \rangle + e_n + \omega \in V_{\frac{1}{n}}(H^{+}({\mathbb{R}}^n))$.\\
The series  $G_{\bf m} (M \langle x \rangle + e_n + \omega_0)$ is bounded on 
$V_{\frac{1}{n}}(H^{+}({\mathbb{R}}^n))$; let us write 
$$
G_{\bf m}(M \langle x \rangle + e_n + \omega_0) \le N.
$$
Since $G_{\bf m}$ is invariant under the action of the translation group 
${\cal{T}}_{n-1}$, the inequality 
$$
\|G_{\bf m}(M \langle x \rangle  + e_n)\| \le N.
$$
thus holds for every $x \in H^{+}(\mathbb{R}^n)$. 
The series $E^{\Lambda}_{s,{\bf m}}(x)$ represent actually  well-defined $s$-monogenic function on the upper half-space. Let us now consider the limit  
\begin{eqnarray*}
\lim\limits_{x_n \rightarrow \infty} E^{\Lambda}_{s,{\bf m}}(x)  &=&   
\sum\limits_{M: T(\Lambda) \backslash \Lambda} \lim\limits_{x_n \rightarrow \infty} 
q^{(s)}_{\bf 0} (c x + d) G_{\bf m}(M \langle x \rangle + e_n)\\
 &=&  \sum\limits_{M: T(\Lambda) \backslash \Lambda,c= 0} \|d\|^{s-n} \lim\limits_{x_n \rightarrow \infty }G_{\bf m}(a[\underline{{\bf x}} + e_n x_n]d^{-1} + \underbrace{bd^{-1}}_{\in \mathbb{Z}^p}\; +e_n)\\
& = & C \lim\limits_{x_n \rightarrow \infty} G_{\bf m}(x)\\
& = & C \sum\limits_{\omega \in \mathbb{Z}^{n-1} \backslash\{0\}} q^{(1)}_{\bf m}(\omega),
\end{eqnarray*}
where $C =\sharp\{M \in \Lambda, c = 0\}$, satisfying $1 \le C \le 2^{p+1}$. The expressions 
\begin{equation}
\label{zeta}
\sum\limits_{\omega \in \mathbb{Z}^{n-1} \backslash\{0\}} q^{(1)}_{\bf m}(\omega), 
\end{equation}
are precisely the Laurent coefficients of the $n-1$-fold periodic monogenic function series 
$$
\epsilon_{\bf m}(z) = \sum\limits_{\omega \in \mathbb{Z}^{n-1}} q^{(1)}_{\bf m}(z+\omega)
$$
which are the partial derivatives of the $n-1$-fold monogenic generalized cotangent function (for details see \cite{KraHabil},~[Chapter 2.1]. There hence must exist multi-indices ${\bf m}$ with $|{\bf m}|\equiv 1(mod\;2)$ for which the associated expression~(\ref{zeta}), 
which can be regarded as a {\it vector-valued generalization of the Riemann zeta function} (see \cite{KraHabil}, [Chapter~2.4]), is different from the zero. Otherwise, one would have $\epsilon_{\bf m}(z) \equiv q^{(1)}_{\bf m}(z)$ which would be a contradiction to the periodicity of $\epsilon_{\bf m}(z)$. For all those indices, the considered limit is hence different from zero, proving the non-triviality of the associated series $E^{\Lambda}_{s,{\bf m}}(x)$ which turns out to be vector-valued by construction. Q.E.D.\\[0.3cm] 
{\bf Odd Weights}\\[0.2cm]
In the cases where $s$ is odd, the construction principle proposed in Theorem~\ref{firstconstruthm} does not produce for all congruence subgroups of $\Gamma_p$ non-trivial examples of automorphic forms. We observe: The groups $\Gamma_p = \Gamma_p[1], \Gamma_p[2]$ as well as all the congruence groups ${\Gamma_p}^{0}[N]$, ${\Gamma_p}_0[N]$ and the generalized theta groups ${\Gamma_p}_{\theta}$ contain the negative identity matrix. Each function that is supposed to satisfy for these groups the transformation behavior   
$$
f(x) = q^{(s)}_{\bf 0}(cx+d) f((ax+b)(cx+d)^{-1}) \quad \quad \mbox{for\;all}\quad M = \left( \begin{array}{cc} a & b \\ c & d \end{array}\right) \in \Lambda
$$
on the whole half-space, has to satisfy in particular 
$$
f(x) = q^{(s)}_{\bf 0}(-x) f((-x)(-x)^{-1}) = -f(x) \quad \quad \mbox{for\;all}\quad x \in H^{+}(\mathbb{R}^n). 
$$
For odd $s$, the proposed construction does hence produce for the groups $\Gamma_p, \Gamma_p[2]$,\\${\Gamma_p}^{0}[N]$, ${\Gamma_p}_0[N], {\Gamma_p}_{\theta}$ only the zero function. 
This effect however does not appear for the smaller principal congruence subgroups $\Gamma_p[N]$ of level $N \ge 3$.\\[0.2cm] 
The next proposition shows that the construction principle proposed in Theorem~\ref{firstconstruthm} provides indeed a number of interesting non-trivial Clifford valued $s$-monogenic automorphic forms for these groups. In what follows ${\cal{T}}_p[N]$ stands for the maximal subgroup of translation matrices that is contained in $\Gamma_p[N]$. 
\begin{proposition}
Suppose that $n,p,N \in {\mathbb{N}}$ with $s<n$ with $n \ge 4$, $p < n-s-1$ and $N \ge 3$. Let $s<n$  Then the following Eisenstein type series 
\begin{equation}
\label{Eisengammapq}
{\cal{G}}^{(p,N)}_{s} (x) = \sum\limits_{M: {\cal{T}}_p[N] \backslash \Gamma_p[N]} q^{(s)}_{\bf 0}(cx+d)
\end{equation} 
represent for all $N \ge 3$ and all $s < n$ (even and odd) non-trivial Clifford-valued monogenic automorphic forms with respect to $\Gamma_p[N]$ on the upper half-space $H^{+}({\mathbb{R}}^n)$. 
\end{proposition}
{\it Proof}. Here again, one can apply the limit argument: 
\begin{eqnarray*}
\lim\limits_{x_n \rightarrow \infty} {\cal{G}}^{(p,N)}_{s} (x)  &=&   \sum\limits_{M: {\cal{T}}_p[n] \backslash \Gamma_p[N]} \lim\limits_{x_n \rightarrow \infty} q^{(s)}_{\bf 0} (c x + d)\\
 &=&  \sum\limits_{M: {\cal{T}}_p[N] \backslash \Gamma_p[N],c = 0} q^{(s)}_{\bf 0}(d) = 1. \hspace{1cm}
\end{eqnarray*}
since $\sharp\{M \in T_p[N] \backslash \Gamma_p[N], c = 0\} = 1$  for all $N \ge 3$. Q.E.D.\\[0.2cm]
To get also non-trivial examples for odd $s$ for the other larger groups, we need to make a more sophisticated construction.  This shall now be explained. One possibility to meet these ends is to consider two weight factors (one from the left and one from the right) and to introduce a second auxiliar variable. Monogenic functions in two vector variables are often called {\it biregular}. For the fundamental theory of biregular functions we refer the reader for example to~\cite{So2}.

From Lemma~1 it readily follows that if a function $f:H^{+}(\mathbb{R}^n) \times H^{+}(\mathbb{R}^n)\rightarrow Cl_n$ satisfies $D_x^s f(x,y) = f(x,y) D_y^t = 0$ for all $(x,y) \in H_2^{+}(\mathbb{R}^n) := H^{+}(\mathbb{R}^n) \times H^{+}(\mathbb{R}^n)$, then 
$$
F(x,y) = q^{(s)}_{\bf 0}(cx+d) f(M\langle x \rangle,M\langle y \rangle) q^{(t)}_{\bf 0}(xc^{*}+d^{*})
$$  
satisfies for all $M \in \Lambda \subseteq \Gamma_p$ the equation $D_x^S F(x,y) = F(x,y) D_y^T = 0$ for all $(x,y) \in H_2^{+}(\mathbb{R}^n) := H^{+}(\mathbb{R}^n) \times H^{+}(\mathbb{R}^n)$ and all $S \ge s, T \ge t$. 
We shall see that it will be advantageous to consider instead the slighly modified transformation  
$$
F(x,y) = \overline{q^{(s)}_{\bf 0}(cx+d)^{*}} f(M\langle x \rangle,M\langle y \rangle) q^{(t)}_{\bf 0}(xc^{*}+d^{*})
$$ 
which has the same analytic properties. 

In order to exclude trivial examples, we however will only regard a function that transforms in the following way 
$$
f(x,y) =  \overline{q^{(s)}_{\bf 0}(cx+d)^{*}} f(M\langle x \rangle,M\langle y \rangle) q^{(t)}_{\bf 0}(yc^{*}+d^{*})
$$
for all $(x,y) \in H_2^{+}(\mathbb{R}^n)$ and all $M \in \Lambda$  
as a {\it non-trivial} modular form on a congruence group $\Lambda\subseteq \Gamma_p$, if additionally its restriction to the diagonal $f(x,x)$ is at least a non-constant well-defined $C^{\infty}$-function that satisfies for all $x \in H^{+}(\mathbb{R}^n)$
$$
\widehat{f}(x):=f(x,x) =  \overline{q^{(s)}_{\bf 0}(cx+d)^{*}} f(M\langle x \rangle,M\langle x \rangle) q^{(t)}_{\bf 0}(xc^{*}+d^{*}) 
$$
for all $M \in \Lambda$. 
The following theorem provides us with a couple of non-trivial examples for all congruence groups $\Lambda \subseteq \Gamma_p$, involving also odd weights. 
\begin{theorem}
\label{secondconstruthm}
Let $p < \min\{n,2n-(s+t)-1\}$ and let $s+t \equiv 0(mod\;2)$. Let $\Lambda \subseteq \Gamma_p$ be an arbitrary congruence saubgroup and $T(\Lambda)$ its subgroup of translation matrices. Suppose that $\tilde{f}: H_2^{+}(\mathbb{R}^n) \rightarrow Cl_n$ is a bounded function satisfying for all $(x,y) \in H_2^{+}(\mathbb{R}^n)$ the equation $D_x^s [\tilde{f}(x,y)] = [\tilde{f}(x,y)] D_y^t = 0$ and additionally $\tilde{f}(T\langle x \rangle, T \langle y \rangle) = T(x,y)$ for all $T \in T(\Lambda)$. Then 
\begin{equation}
f(x,y):= \sum\limits_{M: T(\Lambda) \backslash \Lambda} \overline{q^{(s)}_{\bf 0}(cx+d)^{*}} \tilde{f}(M\langle x \rangle,M\langle y \rangle) q^{(t)}_{\bf 0}(yc^{*}+d^{*})
\end{equation}
is left $s$-monogenic in the variable $x$ and right $t$-monogenic in $y$ and satisfies for all $(x,y) \in H_2^{+}(\mathbb{R}^n)$: 
$$
f(x,y)= \overline{q^{(s)}_{\bf 0}(cx+d)^{*}} f(M\langle x \rangle,M\langle y \rangle) q^{(t)}_{\bf 0}(yc^{*}+d^{*})\quad\quad \mbox{for\;all}\;\;M \in \Lambda.
$$
\end{theorem}
The proof can be done in analogy to that of Theorem~\ref{firstconstruthm} and can be adapted to the setting considered here. In view of the two automorphy factors, we get a better convergence condition:
\begin{eqnarray*}
& & \sum\limits_{M: T(\Lambda) \backslash \Lambda} \|\;\overline{q^{(s)}_{\bf 0}(cx+d)^{*}} \tilde{f}(M \langle x \rangle ,M \langle y \rangle) q^{(t)}_{\bf 0}(y c^{*}+d^{*})\| \\
& \le &  \tilde{L} \rho^{s+t-2n}
\sum\limits_{M:  T(\Lambda) \backslash \Lambda} \frac{1}{\|c e_n + d\|^{n-s}} \frac{1}{\|e_n c^{*} + d^{*}\|^{n-t}}\\
& \le & L  \sum\limits_{M:  T(\Lambda) \backslash \Lambda} \frac{1}{\|c e_n + d\|^{2n-s-t}}.
\end{eqnarray*}
Here $\tilde{L}$ and $L$ denote properly chosen non-negative real constants. 
The series in the previous line thus is absolutely convergent, whenever $p < 2n-s-t-1$. In the monogenic case ($s=t=1$), we hence get convergence for the full hypercomplex modular group $\Gamma_{n-1}$ whenever we are in a space of dimension $n \ge 3$. 
This construction provides us also for odd weights non-trivial polymonogenic automorphic forms on all the congruence subgroups of $\Gamma_p$ where $p<\min\{n,2n-s-t-1\}$; in the monogenic case $s=t=1$ even for all congruence groups of the full modular group $\Gamma_{n-1}$ when we are a space of dimension at least $3$. 
This is established in the following proposition.
\begin{proposition}
Let $n \ge 3$, $s \in \mathbb{N}$ with $s,t < n$.\\
Suppose $p < \min\{n,2n-(s+t)-1\}$ and let $\Lambda \subseteq \Gamma_p$ be an arbitrary congruence subgroup. Then the associated polybiregular Eisenstein series 
\begin{equation}
{\cal{E}}^{\Lambda}_{s}({\bf x},{\bf y}) = \sum\limits_{M: T(\Lambda) \backslash \Lambda} 
\overline{q_{\bf 0}^{(s)}(c_M x+d_M)^{*}} q_{\bf 0}^{(s)}(y c_M^{*}+d_M^{*})
\end{equation} 
are non-trivial Clifford-valued automorphic forms (in the sense of the definition given above) on $\Lambda$ and satisfy on $H_2^{+}(\mathbb{R}^n)$ the relation $D_x^s [{\cal{E}}^{\Lambda}_{s}({\bf x},{\bf y})] = [{\cal{E}}^{\Lambda}_{s}({\bf x},{\bf y})] D_y^s = 0$.  
\end{proposition}
{\it Proof}. To show the non-triviality consider 
\begin{eqnarray}\label{symmetry}
\lim\limits_{x_n \rightarrow \infty} {\cal{E}}_s^{\Lambda}(e_n x_n,e_n x_n) & =&  \sum\limits_{M: T(\Lambda) \backslash \Lambda, c \neq 0} \underbrace{\lim\limits_{x_n \rightarrow \infty} \overline{q^{(s)}_{\bf 0}(c e_n  x_n+d)^{*}}  q^{(s)}_{\bf 0}(e_n x_n c^{*} + d^{*})}_{=0}\nonumber\\
& + & \sum\limits_{M:  T(\Lambda) \backslash \Lambda, c =0} \overline{q^{(s)}_{\bf 0}(d)^{*}}  q^{(s)}_{\bf 0}(d^{*})\nonumber \\
& = & \sum\limits_{M:  T(\Lambda) \backslash \Lambda, c =0} \|q^{(s)}_{\bf 0}(d)\|^2.  
\end{eqnarray}
The last expression equals $\sharp\{M \in T(\Lambda) \backslash \Lambda, c = 0\}$ and is hence a {\it positive} integer smaller or equal than $2^{p+1}$. The upper bound is attained when considering $\Lambda = \Gamma_p$.  

The restriction to the diagonal ${\cal{E}}_s^{\Lambda}(x,x)$ is thus actually a non-constant function. It is furthermore a $C^{\infty}$-function in the single vector variable $x$ and satisfies the desired automorphy relation.\\
Here, the advantage of using $\overline{{q^{(s)}_{\bf 0}}^{*}}$ instead of $q^{(s)}_{\bf 0}$ on the left-hand side becomes clear. The combination of the reversion and the conjugation automorphism on the left automorphy factor causes a symmetry break, so that we get a non-vanishing expression in~(\ref{symmetry}).\\[0.2cm]
The simplest non-trivial Clifford-valued monogenic automorphic forms for an arbitrary congruence subgroup $\Lambda \subseteq \Gamma_{n-1}$, including all principal congruence subgroups $\Gamma_p[N]$ of any level $N \in \mathbb{N}$, the other two basic types of congruence groups ${\Gamma_p}^{0}[N],{\Gamma_p}_{0}[N]$, as well as the generalized theta groups ${\Gamma_p}_{\theta}$ for any $p \le n-1$, are thus given by the following biregular Eisenstein series 
\begin{equation}
{\cal{E}}^{\Lambda}_{1}({\bf x},{\bf y}) = \sum\limits_{M: T(\Lambda) \backslash \Lambda} 
\overline{q_{\bf 0}^{(1)}(c x+d)^{*}} q_{\bf 0}^{(1)}(y c^{*}+d^{*}).
\end{equation} 
{\bf Concluding remarks}: Recent research results indicate that the theory of $s$-monogenic automorphic forms and functions for discrete subgroups of the Vahlen group opens systematically the door to get explicit solutions of a number of boundary value problems on a number of important conformally flat spin manifolds, including for instance the Dirichlet problem.  We refer readers with interest in this branch of applications to the recent paper~\cite{KraRyan2}, jointly written with J.~Ryan, in which an important step in this direction has been done.

Submitted for review: January 9, 2004.

\end{document}